\documentclass[a4paper,12pt]{article}
\usepackage[cp1251]{inputenc}
\usepackage[ukrainian]{babel}
\usepackage{latexsym,amssymb,amsmath,amsthm,amsfonts,euscript}

\newtheorem{theorem}{Теорема}

\newtheorem{lem}[theorem]{Лема}

\newtheorem{definition}{Означення} 

\title{Аналог теореми Меньшова--Трохимчука для моногенних функцій в підпросторі тривимірної комутативної алгебри}

\author{М.В.~Ткачук\\
Інститут математики НАН України\\ 
maxim.v.tkachuk@gmail.com}

\begin{document}

\maketitle

\vskip 5mm

\rightline{UDC 517.54, MSC 30G35}

\begin{abstract}
The aim of this work is to weaken the conditions of monogenity for functions that take values in an arbitrary n-dimensional commutative algebra over the field of complex numbers. The monogenity of the function is understood as a combination of its continuity and the existence of the G\^ateaux derivative.

В статті послаблено умови моногенності функцій зі значеннями в довільній комутативній багатовимірній алгебрі над полем комплексних чисел. Під моногенністю мається на увазі неперервність та існування похідної Гато.
\end{abstract}

{\bf 1. Вступ.} В алгебрі комплексних чисел $\mathbb {C}$ функція
$F\colon \mathbb {C}\longrightarrow\mathbb {C}$ називається
моногенною в точці $\xi_{0}\in \mathbb {C}$, якщо існує скінченна
границя
\begin{equation} \label{compl-der}
\lim _{\xi\to \xi_{0}}{\frac {F(\xi)-F(\xi_{0})}{\xi-\xi_{0}}}\,.
\end{equation}
При цьому границя (\ref{compl-der}) називається похідною функції
$F$ в точці $\xi_{0}$. Функція, яка є моногенною в усіх точках
області $D\subset \mathbb {C}$, називається голоморфною в цій
області (див. \cite{Goursa}).

Встановленню послаблених умов голоморфності функцій комплексної
змінної присвячені роботи Х.~Бора \cite{Bohr}, Х.~Радемахера
\cite{Rademacher}, Д.Є.~Меньшова
\cite{menshov-1,menshov-2,menshov-3}, В.С.~Федорова
\cite{fedorov}, Г.П.~Толстова \cite{Tolstov}, Ю.Ю.~Трохимчука
\cite{Trokhimchuk,zb_trokhinchuk}, Г.Х.~Синдаловського
\cite{Sindalovski,Sindalovski-2},  Д.С.~Теляковського
\cite{Teliakovski, Teliakovski-2, Teliakovski-3, Teliakovski-4},
Є.П.~Долженко \cite{Dolgenko}, М.Т.~Бродович \cite{Brodovich}.

Наведемо тут одну з умов Меньшова, яку, зберігаючи позначення
автора, називають умовою $K'''$, а саме: кажуть, що {\it функція
$F(\xi)$ задовольняє умову $K'''$ в точці $\xi_{0}$, якщо існує
границя (\ref{compl-der}), де $\xi$ належить об'єднанню двох
неколінеарних променів з початком у точці $\xi_{0}$.}

Д.Є.~Меньшов \cite{menshov-1,menshov-2,menshov-3} показав
достатність виконання умови $K'''$ в кожній точці області $D$ (за
винятком не більш ніж зчисленної їх кількості) для конформності
відображення $F$ у випадку, коли $F : D \rightarrow \mathbb{C}$
--- неперервна однолиста функція. Ю.Ю.~Трохимчук
\cite{Trokhimchuk} зняв умову однолистості функції $F$, довівши
при цьому наступну теорему. \vskip 2mm

\noindent {\bf Теорема Меньшова--Трохимчука.} {\it Якщо функція $F
: D \rightarrow \mathbb{C}$ неперервна в області $D$ і в кожній її
точці, за винятком не більш ніж зчисленної їх кількості,
виконується умова $K'''$, то функція $F$ голоморфна в області $D$.
} \vskip 2mm

А.В.~Бондар \cite{bondar,bondar-mono} довів аналог цієї теореми
для функцій, заданих в багатовимірному комплексному просторі
$\mathbb{C}^n$, при цьому ним доведено, що для голоморфності
функції достатньо неперервності цієї функції та існування і
рівності похідної Фреше вздовж $2n$ спеціально вибраних напрямів.
А.В.~Бондар \cite{bondar-mono} і В.І.~Сірик \cite{siryk} довели
також для функцій, заданих в $\mathbb{C}^n$, аналоги іншої теореми
Меньшова--Трохимчука, в якій використовується певна умова
збереження кутів. О.С.~Грецький \cite{gretskii} узагальнив згадані
тут результати А.В.~Бондаря на відображення банахових просторів.

Метою даної роботи є послаблення умов моногенності для функцій, що приймають значення в довільній скінченновимірній комутативній алгебрі над полем комплексних чисел. При цьому моногенність функції розуміється як поєднання її неперервності з існуванням похідної Гато.

\vskip 2mm

{\bf 2. Моногенні функції в багатовимірній комутативній алгебрі.} Розглянемо $n$-вимірну комутативну асоціативну банахову алгебру $\mathbb{A}_n$ з одиницею $1$ над полем $\mathbb{C}$.

Згідно Картана \cite{cartan} існує базис $\{I_k\}_{k=1}^n$ алгебри $\mathbb{A}_n$, який задовольняє наступнив умовам:
\begin{align*}
&1. \qquad I_k I_k = I_k; \; I_k I_s = 0, \; 1 \leq k, s \leq m, k \neq s \\
&2. \qquad I_k I_s = \sum_{j=max\{s,k\}+1}^n \gamma_{k,s}^j I_j, \; m+1 \leq k, s \leq n, \\
\end{align*}

\newpage

базисом якої є трійка $\{1, \rho,
\rho^2\}$, і при цьому виконується рівність $\rho^3 = 0$.
Визначимо евклідову норму елемента алгебри рівністю
$$\|a + b \rho + c \rho^2\|:=\sqrt{|a|^2+|b|^2+|c|^2}\,,\qquad a, b,
c\in \mathbb{C}\,.$$

Алгебра $\mathbb{A}_{3}$ має єдиний максимальний ідеал
$\mathcal{I}:=\{\lambda_{1}\rho+\lambda_{2}\rho^{2}: \lambda_{1},
\lambda_{2}\in\mathbb{C}\}$\, який є також її радикалом.

Оскільки ядром лінійного відображення $f : \mathbb{A}_3
\rightarrow \mathbb{C}$, що визначається рівністю
\begin{equation} \label{mult-func}
f(a + b \rho + c \rho^2) = a\,,
\end{equation}
є максимальний ідеал $\mathcal{I}$, то $f$ є неперервним
мультиплікативним функціоналом (див. \cite[с. 135] {Hil_Filips}).

Зафіксуємо спочатку дійсний тривимірний підпростір  $E_3:=\{\zeta=
xe_1+ye_2+ze_3\,: \,\,x,y,z\in\mathbb{R}\} \subset \mathbb{A}_3$,
де вектори\, $e_1,e_2,e_3$ --- лінійно незалежні над полем дійсних
чисел $\mathbb{R}$, проте, взагалі кажучи, не утворюють базис
алгебри $\mathbb{A}_3$\,. На вибір підпростору $E_3$ накладемо
лише одну вимогу: образом $E_3$ при відображенні $f$ є вся
комплексна площина (див. \cite{Pukh-5,Sh-co}).

Важливими з точки зору застосувань прикладами таких підпросторів є
підпростори, побудовані на гармонічних базисах $\{e_1,e_2,e_3\}$
алгебри $\mathbb{A}_3$\,, що задовольняють рівність\,
$e_1^2+e_2^2+e_3^2=0$ (див.
\cite{melnichenko_plaksa-mono,shpakivskyi_plaksa_umzh}). Існування
гармонічних базисів в комутативній алгебрі є істотною передумовою
побудови розв'язків тривимірного рівняння Лапласа у вигляді
компонент розкладу диференційовних функцій за векторами базису
(див. \cite{Ketchum-28,Mel'nichenko75,melnichenko_plaksa-mono}).

Добре відомо, що існують різні типи диференційовності відображень
в лінійних нормованих просторах.  Насамперед, використовуються
сильна дифекренційовність за Фреше і слабка диференційовність за
Гато  (див., наприклад, \cite{Hil_Filips}), при цьому відповідні
похідні Фреше і Гато визначаються як лінійні оператори.  Для
функції, заданої в області скінченновимірної комутативної
асоціативної алгебри, Г.~Шефферс \cite{Scheffers} розглядав
похідну, яка розуміється як функція, визначена в тій самій
області. Узагальнюючи такий підхід на випадок довільної
комутативної банахової алгебри, Е.Р.~Лорх \cite{Lorch} ввів сильну
похідну функції, яка також розуміється як функція, визначена в тій
же області, що і сама функція.

Функція $\Phi: \Omega \rightarrow \mathbb{A}_3$ називається {\it
диференційовною за Лорхом} в області $\Omega \subset E_3$, якщо
для кожної точки $\zeta \in \Omega$ існує елемент алгебри
$\Phi'_L(\zeta) \in \mathbb{A}_3$ такий, що для кожного
$\varepsilon
> 0$ існує $\delta > 0$ таке, що для всіх $h \in E_3$, для яких
$\|h\| < \delta$, виконується нерівність:
\begin{equation}
\| \Phi(\zeta + h) - \Phi(\zeta) - h \Phi'_L(\zeta) \| \leq \|h\|
\varepsilon.
\end{equation}

Похідна Лорха $\Phi'_L(\zeta)$ є функцією змінної $\zeta$, тобто
$\Phi'_L: \Omega \rightarrow \mathbb{A}_3$.  При цьому
відображення $B_{\zeta}: E_3 \rightarrow \mathbb{A}_3$, задане
рівністю $B_{\zeta} h = h \Phi'_L(\zeta)$, є обмеженим лінійним
оператором. Отже, функція $\Phi$, диференційовна за Лорхом в
області $\Omega$, має похідну Фреше $B_{\zeta}$ в кожній точці
$\zeta \in \Omega$. Обернене твердження загалом не вірне (див.
приклад в монографії \cite[с.~116]{Hil_Filips}).

Використовуючи диференціал Гато,  І.П.~Мельниченко
\cite{Mel'nichenko75} за\-про\-по\-ну\-вав розглядати похідну Гато
також як функцію, визначену в тій же області, що і сама функція.

Якщо для функції $\Phi \colon \Omega\longrightarrow\mathbb{A}_3$,
заданої в області $\Omega\subset E_3$\,, у кожній точці
$\zeta\in\Omega$ існує елемент алгебри
$\Phi_G'(\zeta)\in\mathbb{A}_3$ такий, що
\begin{equation}\label{Gprz}
\lim\limits_{\delta\rightarrow 0+0} \left(\Phi(\zeta+\delta
h)-\Phi(\zeta)\right)\delta^{-1}= h\Phi_G'(\zeta)\quad\forall h\in
E_{3}\,,
\end{equation}
то функцію $\Phi_G' \colon \Omega\longrightarrow\mathbb{A}_3$
будемо називати {\em похідною Гато} функції $\Phi$\,.

Очевидно, що з існування сильної похідної Лорха $\Phi'_L(\zeta)$
випливає існування слабкої похідної Гато $\Phi_G'(\zeta)$ і
рівність $\Phi_L'(\zeta)=\Phi_G'(\zeta)$, проте з існування
похідної Фреше $B_{\zeta}$ не випливає існування похідної
$\Phi_G'(\zeta)$, що демонструє згаданий вище приклад з монографії
\cite[с.~116]{Hil_Filips}.

Розглянемо тепер поняття моногенної функції.

Функцію $\Phi \colon \Omega \longrightarrow \mathbb{A}_3$
називаємо \textit{моногенною} в області $\Omega\subset E_{3}$\,,
якщо $\Phi$ є неперервною і має похідну Гато в кожній точці
області $\Omega$ (див.
\cite{shpakivskyi_plaksa_umzh,overview_plaksa,Plaksa_UMB}).

Хоча з існування похідної Гато $\Phi_G'(\zeta)$ не випливає
існування похідної Лорха $\Phi'_L(\zeta)$\,, але моногенні функції
$\Phi \colon \Omega \longrightarrow \mathbb{A}_3$ в області
$\Omega\subset E_{3}$ є диференційовними за Лорхом у цій області.
Це випливає з представлення моногенних функцій $\Phi(\zeta)$,
$\zeta\in\Omega$\,, через голоморфні функції комплексної змінної
$f(\zeta)$, встановленого в роботі \cite{shpakivskyi_plaksa_umzh}.

В роботі \cite{Pl-zb17} послабено одну з умов моногенності, а
саме: показано, що за умови існування похідної Гато у функції
$\Phi:\Omega\longrightarrow\mathbb{A}_{3}$ в усіх точках області
$\Omega\subset E_{3}$ неперервність функції $\Phi$ можна замінити
її локальною обмеженістю в області $\Omega$\,.

\vskip 2mm

{\bf 3. Аналог теореми Меньшова--Трохимчука для моногенних функцій
в областях простору $E_3$.} 
Введемо деякі
позначення. Перетином радикалу алгебри $\mathbb{A}_3$ з лінійним
простором $E_3$ є множина необоротних елементів, що належать
$E_3$. Цією множиною є деяка пряма $L:=\{c\, l\, :\, c\in
\mathbb{R}\}$, де через $l\in E_3$ позначено напрямний вектор
прямої $L$. Прообразом довільної точки $\xi \in \mathbb{C}$ в
$E_3$ при відображенні $f$ є пряма $L^{\zeta}:=\{\zeta+c\, l\, :\,
c\in \mathbb{R}\}$, де $\zeta$
--- деякий елемент із $E_3$ такий, що $\xi=f(\zeta)$. Очевидно, що
пряма $L^{\zeta}$ паралельна прямій $L$.

Зазначимо, що тут i надалі до об’єктів з $E_3$ застосовуються
геометричні поняття (паралельність, опуклість в напрямку прямої
тощо), які, строго кажучи, мають сенс по відношенню до
конгруентних прообразів цих об’єктів у $\mathbb{R}^3$ при взаємно
однозначній відповідності $\zeta=xe_1+ye_2+ze_3$ між елементами
$\zeta\in E_3$ і точками $(x,y,z)\in\mathbb{R}^3$.

Нехай область $\Omega\subset E_3$ є опуклою в напрямку прямої $L$
(область називається {\it опуклою в напрямку прямої}, якщо вона
містить кожен відрізок, який з’єднує дві точки області i
паралельний цій прямій). При цьому перетини області $\Omega$ з
усіма прямими $L^{\zeta}$, де $\zeta\in\Omega$\,, є зв'язними
внаслідок опуклості області $\Omega$ в напрямку прямої $L$

Розглянемо наступний гіперкомплексний аналог умови Меньшова $K'''$
в алгебрі $\mathbb{A}_3$ для функцій $\Phi : \Omega \rightarrow
\mathbb{A}_3$, визначених в області $\Omega\subset E_3$.

\begin{definition}
Будемо говорити, що функція $\Phi : \Omega \rightarrow
\mathbb{A}_3$ задовольняє умову $K'''_{\mathbb{A}_3, E_3}$ в точці
$\zeta \in \Omega$, якщо існує елемент $\Phi_*(\zeta)\in
\mathbb{A}_3$ такий, що рівність
\begin{equation}
\lim_{\delta \rightarrow 0 + 0} \left( \Phi(\zeta + \delta h) -
\Phi(\zeta) \right)
\delta^{-1} = h \Phi_*(\zeta) 
\label{eq:deriv}
\end{equation}
виконується для трьох векторів $h$, а саме: векторів $h_1, h_2$ і
$h_3 = l$ або $h_3 = -l$, що утворюють базис в просторі $E_3$.
\end{definition}

Зауважимо, що у випадку, коли функція $\Phi : \Omega \rightarrow
\mathbb{A}_3$ задовольняє умову $K'''_{\mathbb{A}_3, E_3}$ в
різних точках області $\Omega \subset E_3$, набір векторів $h_1,
h_2, h_3$ може бути різним в різних точках цієї області.

\begin{lem}
\label{monom_K} Нехай область $\Omega \subset E_3$ є опуклою в
напрямку прямої $L$ і неперервна в $\Omega$ функція $\Phi : \Omega
\rightarrow \mathbb{A}_3$ має вигляд $\Phi(\zeta) = \rho^2
\Phi_2(\zeta)$, де $\Phi_2(\zeta) \in \mathbb{C}$, і задовольняє
умову $K'''_{\mathbb{A}_3, E_3}$ в усіх точках $\zeta \in \Omega$,
крім не більш ніж зчисленної множини точок. Тоді $\Phi_2(\zeta) =
F_2(f(\zeta))$, де $F_2 : D \rightarrow \mathbb{C}$ --- голоморфна
функція в області $D$, 
яка є образом області $\Omega$ при відображенні $f$.
\end{lem}

\begin{proof} Нехай $\zeta \in \Omega$ --- довільна точка, в якій функція $\Phi$ задовольняє умову $K'''_{\mathbb{A}_3,
E_3}$. Запишемо рівність (\ref{eq:deriv}) для функції $\Phi(\zeta)
= \rho^2 \Phi_2(\zeta)$:
\begin{equation}
\lim_{\delta \rightarrow 0 + 0} \rho^2 \left( \Phi_2(\zeta +
\delta h) - \Phi_2(\zeta) \right) \delta^{-1} = h \Phi_*(\zeta)
\label{eq:deriv_monom}
\end{equation}
і зазначимо, що вона виконується при $h\in\{h_1, h_2, h_3\}$.

Підставимо $h = h_1$ у рівність (\ref{eq:deriv_monom}) і з
урахуванням того, що $h_1$ є оборотним елементом алгебри
$\mathbb{A}_3$, отримаємо
\begin{equation} \label{eq:deriv_monom_1}
\Phi_*(\zeta) = \rho^2\, h_1^{-1} \lim_{\delta \rightarrow 0 + 0}
\left( \Phi_2(\zeta + \delta h_1) - \Phi_2(\zeta) \right)
\delta^{-1} =: \rho^2\, \Psi(\zeta).
\end{equation}

Після підстановки виразу (\ref{eq:deriv_monom_1}) для $\Phi_*$ в
рівність (\ref{eq:deriv_monom}) вона набуде вигляду
\begin{equation}
\lim_{\delta \rightarrow 0 + 0} \rho^2 \left( \Phi_2(\zeta +
\delta h) - \Phi_2(\zeta) \right) \delta^{-1} = h \rho^2
\Psi(\zeta). \label{eq:deriv_monom_2}
\end{equation}

Тепер після підстановки в (\ref{eq:deriv_monom_2}) значення $h =
h_3$ отримаємо нуль в правій частині рівності
(\ref{eq:deriv_monom_2}), оскільки $h_3\in\mathcal{I}$\,. Звідси
випливає, що звуження функції $\Phi_2$ на перетин області $\Omega$
з прямою $L^{\zeta}$ в усіх точках, крім не більш ніж зчисленної
множини точок цього перетину, має рівну нулю одну з односторонніх
(взагалі кажучи, різних у різних точках) похідних вздовж прямої
$L^{\zeta}$. При цьому перетин області $\Omega$ з прямою
$L^{\zeta}$ є зв'язним внаслідок опуклості області $\Omega$ в
напрямку прямої $L$. Тоді за теоремою 9 з монографії
Ю.Ю.~Трохимчука \cite[с.103]{zb_trokhinchuk} функція $\Phi_2$ є
сталою на перетині області $\Omega$ з прямою $L^{\zeta}$.

Звідси випливає, що функція $\Phi_2$ може бути представлена у
вигляді $\Phi_2(\zeta) = F_2(f(\zeta))$, де $F_2: D \rightarrow
\mathbb{C}$ --- деяка неперервна в області $D$ функція.

Доведемо, що функція $F_2$ голоморфна в області $D$.

Спочатку зазначимо, що наслідком означення (\ref{mult-func})
функціонала $f$ є рівність
$$ \rho^2\, h \Psi(\zeta)=\rho^2\, f(h) f(\Psi(\zeta)).$$
Тому, позначаючи при цьому $\xi:=f(\zeta)$, переписуємо рівність
(\ref{eq:deriv_monom_2}) у вигляді
\begin{equation}
\rho^2\,\lim_{\delta \rightarrow 0 + 0}
 \left( F_2(\xi + \delta f(h)) - F_2(\xi) \right)
\delta^{-1} = \rho^2\, f(h) f(\Psi(\zeta)).
\label{eq:deriv_monom_3}
\end{equation}

Оскільки вирази біля $\rho^2$ в обох частинах рівності
(\ref{eq:deriv_monom_3}) приймають комплексні значення, то з
єдиності розкладу елемента алгебри за базисом випливає рівність
$$
\lim_{\delta \rightarrow 0 + 0} \left( F_2(\xi + \delta f(h)) -
F_2(\xi) \right) \delta^{-1} = f(h) f(\Psi(\zeta)),
$$
яка виконується при $h \in \{h_1, h_2\}$.

Звідси випливають рівності
\begin{multline*}
f(\Psi(\zeta)) = \lim_{\delta \rightarrow 0 + 0} \left( F_2(\xi +
\delta t_1) - F_2(\xi) \right)
(\delta t_1)^{-1} =\\
=\lim_{\delta \rightarrow 0 + 0} \left( F_2(\xi + \delta t_2) -
F_2(\xi) \right) (\delta t_2)^{-1},
\end{multline*}
де\,\, $t_1 := f(h_1)$, $t_2 := f(h_2)$.

Отже, в кожній точці $\xi$ області $D$, за винятком не більш ніж
зчисленної їх кількості, існують похідні функції $F_2$ вздовж двох
неколінеарних променів з початком у точці $\xi$ і ці похідні
рівні, а це означає, що неперервна функція $F_2$ задовольняє умову
Меньшова $K'''$ у точці $\xi$.  Тоді з теореми
Меньшова--Трохимчука випливає голоморфність функції $F_2$ в
області $D$.
\end{proof}

Кожен елемент $a + b \rho + c \rho^2$, $a, b, c\in \mathbb{C}$, за
умови $a\ne 0$ має обернений елемент, розклад якого за 
базисом $\{1,\rho,\rho^2\}$ визначається рівністю
$$(a + b \rho + c \rho^2)^{-1}=\frac{1}{a}-\frac{b}{a^2}\,\rho +\left(\frac{b^2}{a^3}-\frac{c}{a^2}\right)\rho^2\,.$$
Тоді
\begin{equation}\label{obern-el}
(t - a - b \rho - c \rho^2)^{-1}=\frac{1}{t-a}+\frac{b}{(t-a)^2}\,\rho +\left(\frac{c}{(t-a)^2}+\frac{b^2}{(t-a)^3}\right)\rho^2\,.
\end{equation}
Використовуючи цей розклад, легко виписати розклад за 
базисом $\{1,\rho,\rho^2\}$ головного продовження голоморфної
функції $F : D\rightarrow\mathbb{C}$ в область $\Pi:=\{\zeta\in
E_3 : f(\zeta)\in D\}$, яка очевидно є нескінченним циліндром,
твірні якого паралельні прямій $L$:
\begin{multline}\label{gol-prod}
\frac{1}{2\pi i}\int\limits_{\gamma}F(t)(t-\zeta)^{-1}\,dt=
F(f(\zeta)) +(b_1x+b_2y+b_3z)F'(f(\zeta))\,\rho+\\
+\biggl((c_1x+c_2y+c_3z)F'(f(\zeta))+\frac{(b_1x+b_2y+b_3z)^2}{2}F''(f(\zeta))\biggr)\,\rho^{2}\\
\forall\,\zeta=xe_{1}+ye_{2}+ze_{3}\in\Pi\,,
\end{multline}
 де $i$ --- уявна комплексна одиниця, замкнена жорданова спрямлювана крива $\gamma$
лежить в області $D$ і охоплює точку $f(\zeta)=a_1x+a_2y+a_3z$, а
комплексні сталі $a_k, b_k, c_k$ при $k=1,2,3$ --- це коефіцієнти
з розкладів елементів $e_{1}, e_{2}, e_{3}$ за базисом
$\{1,\rho,\rho^2\}$:
$$\begin{array}{l}
    e_1=a_1+b_1\rho+c_1\rho^2,\\
    e_2=a_2+b_2\rho+c_2\rho^2,\\
    e_3=a_3+b_3\rho+c_3\rho^2.
   \end{array} $$
Розклад (\ref{gol-prod}) узагальнює аналогічний розклад, отриманий
в теоремі 1.7 з \cite{melnichenko_plaksa-mono} при додатковому
припущенні, що $e_1=1$.

\begin{lem}
\label{monom_K-1} Нехай область $\Omega \subset E_3$ є опуклою в
напрямку прямої $L$, функція $\Phi : \Omega \rightarrow
\mathbb{A}_3$ є неперервною в $\Omega$ і задовольняє умову
$K'''_{\mathbb{A}_3, E_3}$ в усіх точках $\zeta \in \Omega$\,,
крім не більш ніж зчисленної множини точок. Тоді при всіх $\zeta
\in \Omega$ справедливе представлення
\begin{equation}\label{predst}
\Phi(\zeta) = \frac{1}{2\pi i}\int\limits_{\gamma}
\Big(F_0(\xi)+F_{1}(\xi)\rho+F_{2}(\xi)\,\rho^2\Big) (\xi -
\zeta)^{-1}\, d \xi\,,
\end{equation}
де $F_0, F_1, F_2$ --- деякі функції, голоморфні в області $D$,
яка є образом області $\Omega$ при відображенні $f$.
\end{lem}

\begin{proof} При $\zeta \in \Omega$ розглянемо розклад $\Phi(\zeta)$ за базисом $\{1,\rho,\rho^2\}$:
$$\Phi(\zeta) = \Phi_0(\zeta) + \Phi_1(\zeta) \rho + \Phi_2(\zeta) \rho^2.$$
 Функція $\rho^2 \Phi(\zeta) = \rho^2
\Phi_0(\zeta)$ є неперервною в $\Omega$ і задовольняє умову
$K'''_{\mathbb{A}_3, E_3}$ в усіх точках $\zeta \in \Omega$, крім
не більш ніж зчисленної множини точок. Тоді з леми \ref{monom_K}
випливає, що $\Phi_0(\zeta) = F_0(f(\zeta))$, де $F_0$ ---
голоморфна функція в області $D$, яка є образом області $\Omega$
при відображенні $f$.

Як випливає з рівності (\ref{gol-prod}), перші компоненти в
розкладах за базисом $\{1,\rho,\rho^2\}$ функцій $\Phi(\zeta)$ і
$\frac{1}{2\pi i}\int_{\gamma} F_0(\xi) (\xi - \zeta)^{-1}\, d
\xi$ співпадають в області $\Omega$. Тому справедлива рівність
\begin{equation}\label{rivn1}
\Phi(\zeta) - \frac{1}{2\pi i}\int\limits_{\gamma} F_0(\xi) (\xi -
\zeta)^{-1}\, d \xi = \Phi_{11}(\zeta)\, \rho + \Phi_{12}(\zeta)\,
\rho^2 \quad \forall\zeta \in \Omega,
\end{equation}
 де
$\Phi_{11}$, $\Phi_{12}$ --- деякі комплекснозначні неперервні в
$\Omega$ функції.

Тоді функція $\rho (\Phi_{11}(\zeta) \rho + \Phi_{12}(\zeta)
\rho^2) = \rho^2 \Phi_{11}(\zeta)$  є неперервною в $\Omega$ і
задовольняє умову $K'''_{\mathbb{A}_3, E_3}$ в усіх точках $\zeta
\in \Omega$, крім не більш ніж зчисленної множини точок. Отже, за
лемою \ref{monom_K} маємо $\Phi_{11}(\zeta) = F_1(f(\zeta))$, де
$F_1$
--- голоморфна функція в області $D$.

Далі так, як і при доведенні рівності (\ref{rivn1}), отримуємо
рівність
\begin{equation}\label{rivn2}
\Phi_{11}(\zeta)\, \rho + \Phi_{12}(\zeta)\, \rho^2 -
\rho\,\frac{1}{2\pi i}\int\limits_{\gamma} F_1(\xi) (\xi -
\zeta)^{-1}\, d \xi = \Phi_{22}(\zeta)\, \rho^2 \quad \forall\zeta
\in \Omega,
\end{equation}
 де $\Phi_{22}$ --- деяка комплекснозначна неперервна в
$\Omega$ функція.

Як наслідок рівностей (\ref{rivn1}), (\ref{rivn2}), маємо рівність
\begin{multline}\label{rivn3}
\Phi(\zeta) - \frac{1}{2\pi i}\int\limits_{\gamma} F_0(\xi) (\xi -
\zeta)^{-1}\, d \xi -\\
-\rho\, \frac{1}{2\pi
i}\int\limits_{\gamma} F_{1}(\xi) (\xi - \zeta)^{-1}\, d \xi =
\Phi_{22}(\zeta) \rho^2 \quad \forall\zeta \in \Omega.
\end{multline}

Тепер, спираючись на лему \ref{monom_K}, приходимо до рівності
$\Phi_{22}(\zeta) = F_2(f(\zeta))$, де $F_2$
--- голоморфна функція в області $D$. Тому справедливими є також
рівності
\begin{equation}\label{rivn4}
\rho^2\,\Phi_{22}(\zeta)=\rho^2\,F_2(f(\zeta))  = \rho^2\,
\frac{1}{2\pi i}\int\limits_{\gamma} F_{2}(\xi) (\xi -
\zeta)^{-1}\, d\xi \quad \forall\zeta \in \Omega\,.
\end{equation}

Нарешті, як наслідок рівностей (\ref{rivn3}), (\ref{rivn4}),
отримуємо представлення (\ref{predst}).
\end{proof}

Основним результатом пункту 3 є наступне твердження.

\begin{theorem} \label{theor-1}
Нехай область $\Omega \subset E_3$ є опуклою в напрямку прямої
$L$, функція $\Phi : \Omega \rightarrow \mathbb{A}_3$ є
неперервною в $\Omega$ і задовольняє умову $K'''_{\mathbb{A}_3,
E_3}$ в усіх точках $\zeta \in \Omega$\,, крім не більш ніж
зчисленної множини точок. Тоді:

1) функція $\Phi$ є моногенною в області $\Omega$;

2) функція $\Phi$ продовжується до функції, моногенної в області
$\Pi$. Таке продовження єдине і задається рівністю (\ref{predst})
при всіх $\zeta\in\Pi$;

3) моногенне продовження (\ref{predst}) функції $\Phi$  є
диференційовним за Лорхом в області $\Pi$.
\end{theorem}

Усі твердження теореми \ref{theor-1} є очевидними наслідками
пред\-став\-лен\-ня (\ref{predst}).

\vskip 2mm

{\bf 4. Узагальнення на інші розмірності області визначення
функцій.} Отримані результати узагальнимо на дійсний підпростір
$E_k$ алгебри $\mathbb{A}_3$ довільної розмірності $2 \le k \leq
6$, на вибір якого накладемо лише одну вимогу: образом $E_k$ при
відображенні $f$ є вся комплексна площина.

Для цього сформулюємо аналог умови Меньшова $K'''$ для функцій
$\Phi : \Omega \rightarrow \mathbb{A}_3$, визначених в області
$\Omega\subset E_k$. Зазначимо, що перетином простору $E_k$ з
радикалом $\mathcal{I}$ алгебри $\mathbb{A}_3$ є площина
розмірності $k-2$\,, яку позначимо $L_{E_k}$. Зокрема, $L_{E_3}$
--- це пряма $L$, визначена в пункті 3, і\, $L_{E_2}=\{0\}$\,.

\theoremstyle{definition}
\begin{definition}
Будемо говорити, що функція $\Phi : \Omega \rightarrow
\mathbb{A}_3$ задовольняє умову $K'''_{\mathbb{A}_3, E_k}$ в точці
$\zeta \in \Omega\subset E_k$, якщо існує елемент
$\Phi_*(\zeta)\in \mathbb{A}_3$ такий, що рівність
(\ref{eq:deriv}) виконується для $k$ різних векторів, а саме: двох
векторів $h_1, h_2$, що мають неколінеарні образи при відображенні
$f$, і $k-2$ векторів $h_3, ... , h_k$, що утворюють базис в
просторі $L_{E_k}$.
\end{definition}

Зауважимо, що у випадку, коли функція $\Phi : \Omega \rightarrow
\mathbb{A}_3$ задовольняє умову $K'''_{\mathbb{A}_3, E_k}$ в
різних точках області $\Omega \subset E_k$, набір векторів $h_1,
h_2,\dots, h_k$ може бути різним в різних точках цієї області.

Наступне твердження є узагальненням леми \ref{monom_K}.

\begin{lem} \label{gen-lem1}
Нехай область $\Omega \subset E_k$ має зв'язні перетини площинами
$L_{E_k}^{\zeta}:=\{\zeta+\tau\, :\, \tau\in L_{E_k}\}$, де
$\zeta\in\Omega$, паралельними площині $L_{E_k}$, і неперервна в
$\Omega$ функція $\Phi : \Omega \rightarrow \mathbb{A}_3$ має
вигляд $\Phi(\zeta) = \rho^2 \Phi_2(\zeta)$, де $\Phi_2(\zeta) \in
\mathbb{C}$, і задовольняє умову $K'''_{\mathbb{A}_3, E_k}$ в усіх
точках $\zeta \in \Omega$, крім не більш ніж зчисленної множини
точок. Тоді $\Phi_2(\zeta) = F_2(f(\zeta))$, де $F_2 : D
\rightarrow \mathbb{C}$ --- голоморфна функція в області $D$, яка
є образом області $\Omega$ при відображенні $f$.
\end{lem}

\begin{proof} Нехай $\zeta \in \Omega$ --- довільна точка, в якій функція $\Phi$ задовольняє умову $K'''_{\mathbb{A}_3, E_k}$\,,
Так, як і при доведенні леми \ref{monom_K}, отримуємо рівність
(\ref{eq:deriv_monom_2}), яка виконується при $h\in\{h_1,
h_2,\dots, h_k\}$. Тепер при підстановці в рівність
(\ref{eq:deriv_monom_2}) значень $h = h_3, ... , h_k$ отримуємо
нуль в правій частині цієї рівності.

Звідси випливає, що похідні функції $\Phi_2$ вздовж напрямків
$h_3, ... , h_k$ (які, взагалі кажучи, є різними в різних точках)
рівні нулю всюди на множині $L_{E_k}^{\zeta}\cap\Omega$\,, крім не
більш ніж зчисленної множини точок. При цьому множина\,
$L_{E_k}^{\zeta}\cap\Omega$\, є зв'язною. Тоді за теоремою 9 з
монографії Ю.Ю.~Трохимчука \cite[с.103]{zb_trokhinchuk} функція
$\Phi_2$ є сталою на множині $L_{E_k}^{\zeta}\cap\Omega$\,.

Звідси випливає, що функція $\Phi_2$ може бути представлена у
вигляді $\Phi_2(\zeta) = F_2(f(\zeta))$, де $F_2: D \rightarrow
\mathbb{C}$ --- деяка неперервна в області $D$ функція.

Доведення голоморфності функції $F_2$ в області $D$ аналогічне
доведенню голоморфності функції $F_2$ в лемі \ref{monom_K}.
\end{proof}

Зауважимо, що у випадку $k=2$ в лемі \ref{gen-lem1} виконується
рівність $L_{E_2}^{\zeta}=\{\zeta\}$ і тому умова про зв'язність
перетинів $L_{E_2}^{\zeta}\cap\Omega$ при $\zeta \in \Omega$,
очевидно, виконується автоматично.

Наступне твердження є узагальненням теореми \ref{theor-1} на
випадок функцій, визначених в областях дійсного підпростору $E_k$
алгебри $\mathbb{A}_3$ довільної розмірності $2 \le k \leq 6$.

\begin{theorem} \label{theor-2}
Нехай область $\Omega \subset E_k$ має зв'язні перетини площинами
$L_{E_k}^{\zeta}$, де $\zeta\in\Omega$, паралельними площині
$L_{E_k}$, і неперервна в $\Omega$ функція $\Phi : \Omega
\rightarrow \mathbb{A}_3$ задовольняє умову $K'''_{\mathbb{A}_3,
E_k}$ в усіх точках $\zeta \in \Omega$\,, крім не більш ніж
зчисленної множини точок. Тоді:

1) функція $\Phi$ є моногенною в області $\Omega$;

2) функція $\Phi$ продовжується до функції, моногенної в області
$\Pi:=\{\zeta\in E_k : f(\zeta)\in D\}$. Таке продовження єдине і
задається рівністю (\ref{predst}) при всіх $\zeta\in\Pi$;

3) моногенне продовження (\ref{predst}) функції $\Phi$  є
диференційовним за Лорхом в області $\Pi$.
\end{theorem}

\begin{proof}
Використовуючи рівність (\ref{obern-el}), отримуємо розклад за
ба\-зисом $\{1,\rho,\rho^2\}$ головного продовження довільної
голоморфної функції $F : D\rightarrow\mathbb{C}$ в область $\Pi$:
\begin{multline}\label{gol-prod-k}
\frac{1}{2\pi
i}\int\limits_{\gamma}F(t)(t-\zeta)^{-1}\,dt= F(f(\zeta)) 
+(b_1x_1+b_2x_2+...+b_kx_k)F'(f(\zeta))\,\rho+\\
+\biggl((c_1x_1+c_2x_2+...+c_kx_k)F'(f(\zeta))+\frac{(b_1x_1+b_2x_2+...+b_kx_k)^2}{2}F''(f(\zeta))\biggr)\,\rho^{2}\\
\forall\,\zeta=x_1e_{1}+x_2e_{2}+...+x_ke_{k}\in\Pi\,,
\end{multline}
 де
замкнена жорданова спрямлювана крива $\gamma$ лежить в області $D$
і охоплює точку $f(\zeta)=a_1x_1+a_2x_2+...+a_kx_k$, а комплексні
сталі $a_j, b_j, c_j$ при $j=1,2,...,k$ --- це коефіцієнти з
розкладів елементів $e_{1}, e_{2}, ..., e_{k}$ за базисом
$\{1,\rho,\rho^2\}$:
$$\begin{array}{l}
    e_1=a_1+b_1\rho+c_1\rho^2,\\
    e_2=a_2+b_2\rho+c_2\rho^2,\\
    ................................\\
    e_k=a_k+b_k\rho+c_k\rho^2.
   \end{array} $$

Тоді, спираючись на лему \ref{gen-lem1}, так, як і при доведенні
леми \ref{monom_K-1}, отримуємо інтегральне представлення
(\ref{predst}) функції $\Phi$, з якого випливають усі твердження
теореми \ref{theor-2}.
\end{proof}

Частина результатів роботи анонсована в препринті
\cite{Tk-Pl-prep}.

\end{document}